\documentclass[12pt]{article}
\usepackage[english]{style}
\usepackage{latexsym,amsmath,amssymb,amsfonts,graphics,xypic}

\newcommand{\shrinkmargins}[1]{
  \addtolength{\textheight}{#1\topmargin}
  \addtolength{\textheight}{#1\topmargin}
  \addtolength{\textwidth}{#1\oddsidemargin}
  \addtolength{\textwidth}{#1\evensidemargin}
  \addtolength{\topmargin}{-#1\topmargin}
  \addtolength{\oddsidemargin}{-#1\oddsidemargin}
  \addtolength{\evensidemargin}{-#1\evensidemargin}
  }

\shrinkmargins{.7}

\DeclareFontFamily{U}{rsf}{}
\DeclareFontShape{U}{rsf}{m}{n}{
  <5> <6> rsfs5 <7> <8> <9> rsfs7 <10-> rsfs10}{}
\DeclareMathAlphabet{\mathscr}{U}{rsf}{m}{n}

\DeclareMathAlphabet{\mathgth}{U}{euf}{m}{n}

\DeclareFontFamily{U}{cyr}{}
\DeclareFontShape{U}{cyr}{m}{n}{
  <5> wncyr5 <6> wncyr6 <7> wncyr7 <8> wncyr8 <9> wncyr9 <10-> wncyr10}{}
\DeclareMathAlphabet{\mathcyr}{U}{cyr}{m}{n}

\input cyracc.def

\makeatletter
\def\operator@font{\sf}
\makeatother

\newcommand{\sL}{{\mathcal L}}

\newcommand{\cJ}{{\mathscr J}}

\newcommand{\cO}{{\mathscr O}}

\newcommand{\R}{\mathbf{R}}

\renewcommand{\P}{{{\mathbf{P}}}} 
\newcommand{\field}[1]{\mathbb{#1}}
\newcommand{\F}{\field{F}}
\renewcommand{\phi}{\varphi}

\title{On the joints problem with multiplicities}
\author{M\'arton Hablicsek}
\begin{document}

\maketitle

\section{Introduction}

\paragraph
Let $k$ be an arbitrary field and consider a finite set $\sL$ of lines in $k^n$  ($n\geq 3$). We denote by $L$ the size of $\sL$. We say that $x\in k^n$ is a joint formed by $\sL$ if there exist at least $n$ lines going through $x$ whose direction vectors are linearly independent. We denote the set of joints by $\cJ$. For a joint $x$ we denote by $N(x)$ the number of $n$-tuples of lines in $\sL$ which form a joint at $x$; and by $r(x)$ the number of lines in $\sL$ passing though $x$.

A simple construction shows that the number of joints, $J$, can be as big as $O_n(L^{\frac{n}{n-1}})$. In a groundbreaking paper, Guth and Katz (\cite{GuKa}) proved that this is indeed the upper bound for the case $n=3$, $k=\R$. Recently Kaplan, Sharir and Shustin \cite{KaShSh}, Quilodr\'an \cite{Qu}, Dvir \cite{Dv}, Tao \cite{Tao} and Carbery and Iliopoulou \cite{CaIl} simplified the proof of Guth and Katz and extended the result to any field $k$ and to any dimension $n$:

\begin{Theorem}[Joints problem]
In $k^n$ we have
\[J=O_n(L^{\frac{n}{n-1}}).\]
\end{Theorem}
\medskip

\noindent
At IPAM in 2014, Carbery and Iliopoulou asked for a simple proof using the polynomial method of the following variant of the joints problem.

\begin{Conjecture}
\label{con:main}
Let $\sL$ be a finite set of lines in $k^n$ of size $L$. Then, the number of joints counted with multiplicities satisfies
\[\sum_{x\in \cJ}N(x)^{\frac{1}{n-1}}\leq c_nL^{\frac{n}{n-1}}\]
where $c_n$ is a constant depending only on $n$.
\end{Conjecture}
\medskip

The purpose of this short note is to give a proof of Conjecture \ref{con:main}, assuming an extra hypothesis: at each joint $x$ any $n$-tuple of lines in $\sL$ passing through $x$ form a joint (in particular, we have $N(x)={r(x) \choose n}$):

\begin{Theorem}
\label{thm:main}
Let $\sL$ be a finite set of lines in $k^n$ of size $L$. Assume that at each joint $x$ any $n$-tuple of lines in $\sL$ passing through $x$ form a joint. Then,
\[\sum_{x\in \cJ}N(x)^{\frac{1}{n-1}}\leq c_nL^{\frac{n}{n-1}}\]
where $c_n$ is a constant depending only on $n$.
\end{Theorem}
\medskip

Without this extra hypothesis the above conjecture is still open except in the case of $n=3$ and $k=\R$ which was solved by Iliopoulou (see \cite{Ili12}, \cite{Ili13}). We also remark that a slightly stronger version of Theorem \ref{thm:main} was proved independently by Iliopoulou (\cite{Ili14}) in the $k=\R$ case.

We also remark that without loss of generality we can assume that the field $k$ is algebraically closed: $k$-lines intersect the same way in $k^n$ and in $\bar{k}^n$ (where $\bar{k}$ is the algebraic closure of $k$). In the sequel $k$ denotes an algebraically closed field.

\section{Acknowledgment}
The author grateful to Anthony Carbery for his encouragement and for his useful remarks on an earlier draft. The author also would like to thank Marina Iliopoulou for useful discussions and for sharing her recent work. Part of this research was performed while the author was visiting the Institute for Pure and Applied Mathematics (IPAM), which is supported by the National Science Foundation. 

\section{Preliminaries}

In this section we collect some geometric facts we need later. We begin with two standard lemmas bounding the degree of non-zero polynomials vanishing on finite sets.

\begin{Lemma}
\label{lem:polmet} Let $\cJ$ be a finite set of points in $k^n$. Let ${m(x)}_{x\in \cJ}$ be a finite collection of natural numbers. Then, there exists a non-zero polynomial $p$ of degree at most $(\sum_{x\in \cJ}(m(x)+n)^n)^{\frac{1}{n}}$ vanishing at the points $x\in \cJ$ to order at least $m(x)$.

In particular, if for all $x\in \cJ$, $m(x)\geq n$, then there exists a non-zero polynomial  of degree at most $2(\sum_{x\in \cJ}m(x)^n)^{\frac{1}{n}}$ vanishing at the points $x\in \cJ$ to order at least $m(x)$.
\end{Lemma}
\medskip

\begin{Proof}
The polynomials of degree at most $d$ form a vector space of dimension ${d+n\choose n}$. A polynomial vanishes at a point $x$ to order at least $m(x)$ if all the derivatives\footnote{In positive characteristic one needs to be more careful and should consider the so-called Hasse derivatives (see, for example, Lemma 2.3 of \cite{Tao}).} of the polynomial of order less than $m(x)$ vanish. Therefore, if
\[\sum_{x\in \cJ}{m(x)+n-1\choose n}<{d+n \choose n}\]
we obtain a non-zero polynomial vanishing at all the $x$ of order of vanishing at least $m(x)$. The inequalities 
\[\sum_{x\in \cJ}{m(x)+n-1\choose n}<\sum_{x\in \cJ}\frac{(m(x)+n)^n}{n!}\]
and
\[{d+n\choose n}>\frac{d^n}{n!}\]
imply the statement of the above lemma.\qed
\end{Proof}
\medskip

\begin{Lemma}
\label{lem:lines}
Let $\sL$ be a set of lines in $k^n$ of size $L$. Then, there exists a non-zero polynomial $p$ of degree at most $nL^{\frac{1}{n-1}}$ vanishing on the lines. 
\end{Lemma}
\medskip

\begin{Proof}
The polynomials of degree at most $d$ form a vector space of dimension 
\[{d+n\choose n}> (d+1)\cdot \frac{d^{n-1}}{n!}.\]
Pick $d+1$ points on each line. If the inequality 
\begin{align}
\label{ine:lem2}
(d+1)L<(d+1)\frac{d^{n-1}}{n!}
\end{align}
holds, then the inequality $(d+1)L<{d+n \choose n}$ holds as well showing that there exists a non-zero degree $d$ polynomial vanishing on all the points. Since the polynomial vanishes on $d+1$ points on each line, by B\'ezout's theorem the polynomial vanishes on all the lines. 

\noindent
From inequality \ref{ine:lem2} we obtain $L< \frac{d^{n-1}}{n!}$, which holds for $d=nL^{\frac{1}{n-1}}$.\qed
\end{Proof}
\medskip

\noindent
We conclude this section by showing a slight generalization of Proposition 13 of \cite{Ko}.

\begin{Lemma}
\label{lem:Kol}
Let $S_1,\dots, S_{n-1}$ be hypersurfaces of degrees $a_i$ in the projective space $\P_k^n$. Assume that the hypersurfaces have no common irreducible components. Let $\sL$ be a union of lines contained in the curve $S_1\cap \dots \cap S_{n-1}$. We denote by $r(x)$ the number of lines in $\sL$ passing through $x\in S_1\cap\dots \cap S_{n-1}$. Then, there exists a constant $M$ depending on only $n$, so that
\[\sum_{x} r(x)^{\frac{n}{n-1}}\leq \sum_{i=1}^{n-1}a_i\cdot \prod_{i=1}^{n-1} a_i\]
where the summation is over those $x\in S_1\cap\dots \cap S_{n-1}$ for which $r(x)>M$.
\end{Lemma}
\medskip

\begin{Proof}
This is a direct adaptation of Proposition 13 of \cite{Ko}, for sake of completeness, we show the key steps. The Hilbert polynomial of a complete intersection curve $B=S_1\cap\dots\cap S_{n-1}$ is given by
\[H_B(t)=\prod_i a_i \cdot t-\frac{1}{2}(\sum_i a_i-n-1)\cdot \prod_i a_i.\]
The constant term, $\frac{1}{2}(\sum_i a_i-n-1)\cdot \prod_i a_i$, is related to the arithmetic genus of $B$,
\[p_a(B)=1+\frac{1}{2}(\sum_i a_i-n-1)\cdot \prod_i a_i.\]
For complete intersection curves $B$, the arithmetic genus is equal to the dimension of the first cohomology space of $\cO_B$,
\[p_a(B)=h^1(B,\cO_B).\]
We compare $h^1(B,\cO_B)$ with $h^1(C,\cO_C)$, where $C$ denotes the reduced subcurve which is the union of lines contained in $S_1\cap\dots\cap S_{n-1}$. On one hand, a basic sheaf theoretic argument shows that
\[h^1(C,\cO_C)\leq h^1(B,\cO_B)\]
and thus
\begin{align}
\label{eq:1}
h^1(C,\cO_C)&\leq 1+\frac{1}{2}(\sum_i a_i-n-1)\cdot \prod_i a_i< \frac{1}{2}(\sum_i a_i)\cdot \prod_i a_i.
\end{align}
On the other hand, it is shown in \cite{Ko}, that
\begin{align}
\label{eq:2}
h^1(C,\cO_C)&\geq \sum_x \delta^*(x)
\end{align}
where the summation is over the singular points\footnote{We say that a point on a scheme $X$ is singular it the dimension of the tangent space at that point is bigger than the dimension of $X$.} of $C$, and $\delta^*(x)$ is called the genus of the singularity. The points $x$ for which $r(x)$ is at least 2 are singular points of $C$. It is also shown in \cite{Ko}, that for such points the genus of the singularity satisfies
\[\delta^*(x)\geq\sum_i \left(r(x)-{i+n-1 \choose n-1}\right)\]
where the summation is over those $i\geq 0$ for which $r(x)-{i+n-1 \choose n-1}>0$. If $r(x)$ is large enough, then 
\[\delta^*(x)\geq \frac{1}{2}r(x)^{\frac{n}{n-1}}.\]
The above inequality combined with inequalities \ref{eq:1} and \ref{eq:2} imply the statement of the lemma.\qed
\end{Proof}
\medskip

\noindent
We remark that the above lemma clearly remains true for any set of points $A\subset S_1\cap\dots\cap S_{n-1}$ satisfying $r(x)>M$ for all $x\in A$. 

\section{Proof of the main theorem}

In this section we prove the main theorem which we restate below.

\begin{Theorem}
\label{thm:main}
Let $\sL$ be a finite set of lines in $k^n$ of size $L$ where $k$ is an arbitrary field, and let $\cJ$ be the set of joints formed by $\sL$. For a joint $x$ we denote by $r(x)$ the number of lines in $\sL$ passing through $x$. Assume that at each joint $x$ any $n$ lines of $\sL$ passing through $x$ form a joint at $x$. Then,
\[\sum_{x\in\cJ} r(x)^{\frac{n}{n-1}}\leq C_n L^{\frac{n}{n-1}}\]
where $C_n$ is a constant depending only on $n$.
\end{Theorem}

\begin{Proof}
We proceed by contradiction. Assume that there exists a configuration so that
\[\sum_{x\in \cJ}r(x)^{\frac{n}{n-1}}=K\cdot L^{\frac{n}{n-1}}\]
for some constant $K$ which we choose later. We also fix a large constant $M$. Without loss of generality, we can assume that for every joint $x$, $r(x)> M$. Indeed, let us denote the set of joints where $r(x)\leq M$ by $\cJ_M$. Then,
\[\sum_{x\in \cJ_M}r(x)^{\frac{n}{n-1}}\leq M^{\frac{n}{n-1}}\sum_{x\in \cJ_M}1\leq M^{\frac{n}{n-1}}|\cJ_M|\leq M^{\frac{n}{n-1}}L^{\frac{n}{n-1}}.\]
If $K$ is much larger than $M^{\frac{n}{n-1}}$, then the contribution from $\cJ_M$ is negligible; thus from now on we assume that $r(x)>M$ for every joint $x$.
\medskip

{\em Step} 1:
We choose $M$ to be large enough (at least $n$) to apply the second part of Lemma \ref{lem:polmet}. There exists a non-zero polynomial $q$ of degree at most 
\[d\leq 2(\sum_{x\in \cJ}r(x)^{\frac{n}{n-1}})^{\frac{1}{n}}\]
vanishing on all the joints to order at least $r(x)^{\frac{1}{n-1}}$. The polynomial $q$ may not be irreducible. We write $q$ as a product of its irreducible factors $\prod_i q_i$. We denote the degree of $q_i$ by $d_i$, in particular, we have $d=\sum_i d_i$. 
Let $m_i(x)$ denote the order of vanishing of $q_i$ at the point $x$. We know that 
\[\sum_i m_i(x)\geq r(x)^{\frac{1}{n-1}}.\]
We choose non-negative numbers $n_i(x)$ at each joint, satisfying 
\begin{itemize}
\item $\sum_i n_i(x)=r(x)^{\frac{1}{n-1}}$, and
\item $n_i(x)\leq m_i(x)$.
\end{itemize}
For subsets $\sL'\subset \sL$ and $\cJ'\subset \cJ$ we define the weighted incidence count as \[I_i(\sL',\cJ'):=\sum_{l\in \sL'}\sum_{x\in l\cap \cJ'}n_i(x).\]
It is easy to see that 
\[\sum_i I_i(\sL,\cJ)=\sum_{x\in \cJ}r(x)^{\frac{n}{n-1}}=KL^{\frac{n}{n-1}}.\]
Since $\sum_i d_i=d\leq 2(\sum_{x\in \cJ} r(x)^{\frac{n}{n-1}})^{\frac{1}{n}}$, there exists an $i$ so that
\[I_i(\sL,\cJ)\geq d_i (\sum_{x\in \cJ} r(x)^{\frac{n}{n-1}})^{\frac{n-1}{n}}/2=d_iK^{\frac{n-1}{n}}L/2.\]
We denote the hypersurface corresponding to $q_i$ by $S_1$.
\medskip

{\em Step} 2:
We define $\sL_i$ to be the set of lines incident to at least $3d_i$ points with respect to the weighted incidence count:
\[\sL_i:=\{l\in \sL:I_i(l,\cJ)\geq 3d_i\}.\]
We note that any line in $\sL_i$ is contained in $S_1$. Furthermore, we define $\cJ_i$ to be the set of those points in $\cJ$ which lie on at least $M/3$ lines of $\sL_i$. We note that if $M$ is large enough ($M\geq 3n$) then any such point is a joint with respect to $\sL_i$, in particular it is singular point of $S_1$ (at each such point there are at least $n$ lines contained in $S_1$ spanning $k^n$).

As above, we define $\sL_i'$ to be the set of those lines which are incident to more than $2d_i$ points of $\cJ_i$ with respect to the weighted incidence count:
\[\sL_i':=\{l\in \sL_i:I_i(l,\cJ_i)> 2d_i\}.\]
Finally, we define $\cJ_i'$ to be the set of those joints which lie on at least $M/3$ lines of $\sL_i'$. Notice that if $K$ is large enough, then $I_i(\sL_i',\cJ_i')\geq d_iK^{\frac{n-1}{n}}L/8$. Indeed:
\[I_i(\sL_i',\cJ_i')= I_i(\sL,\cJ)-I_i(\sL\setminus \sL_i,\cJ)-I_i(\sL_i,\cJ\setminus \cJ_i)-I_i(\sL_i\setminus \sL_i',\cJ_i)-I_i(\sL_i',\cJ_i\setminus \cJ_i')\]
where 
\[I_i(\sL\setminus \sL_i,\cJ)+I_i(\sL_i\setminus \sL_i',\cJ_i)\leq 5d_iL\]
by definition; and since $r(x)>M$, we have
\[I_i(\sL_i,\cJ\setminus \cJ_i)\leq \sum_{x\in \cJ\setminus \cJ'}n_i(x)M/3\leq \sum_{x\in \cJ}n_i(x)r(x)/3=I_i(\sL,\cJ)/3;\]
and similarly
\[I_i(\sL_i',\cJ_i\setminus \cJ_i')\leq I_i(\sL,\cJ)/3.\]
Summarizing the above discussion we obtain
\begin{align}
\label{eq:main}
\sum_{x\in \cJ_i'}n_i(x)r_i(x)\geq d_iK^{\frac{n-1}{n}}L/8
\end{align}
where $r_i(x)>M/3$ denotes the number of lines of $\sL_i'$ passing through $x$.
\medskip

{\em Step} 3:
Consider the gradient of $q_i$: $\nabla q_i$. Since every $x\in \cJ_i$ is a singular point of $S_1$, therefore the components of the gradient vanish at those points to order at least 
\[\max(n_i-1,1)\geq n_i/2.\] 
Recall that $I(l,\cJ_i)> 2d_i$ for any $l\in \sL_i'$, hence $\sL_i'$ is contained in the vanishing locus of any component of the gradient. The vanishing on the gradient $\nabla q_i$ implies that $q_i$ is constant\footnote{In characteristic $p>0$ the vanishing of the gradient implies that $q_i=f^p$ for some other polynomial $f$. Since $k$ is algebraically closed, this contradicts the assumption that $q_i$ is irreducible unless $f$ is constant.}. On the other hand the polynomial vanishes on the joints $x\in \cJ'$ to order at least $n_i(x)/2$, hence  we can pick a component of its gradient which does not vanish on $S_1$: we denote it by $q_i'$.
\medskip 

{\em Recursive step}:
We note that the above set up is exactly the same set up we started with: 
\begin{itemize}
\item $\sL_i'$ plays the role of $\sL$,
\item $\cJ_i'$ plays the role of $\cJ$,
\item $q_i'$ plays the role of $q$, 
\item $d_i$ plays the role of $d$,
\item $M/3$ plays the role of $M$, 
\item $n_i(x)/2$ plays the role of $m_i(x)$,
\item $r_i(x)$ plays the role of $r(x)$ (it satisfies $r_i(x)>M/3$ for all $x\in\cJ_i'$ similarly to the inequality $r(x)>M$ satisfied for all $x\in \cJ$).
\end{itemize}
We proceed as follows. We write $q_i'$ as a product of irreducible factors $q_{ij}$. We denote the degree of $q_{ij}$ by $d_{ij}$; in particular, we have
\[\sum_j d_{ij}\leq d_i\leq d\leq 2K^{\frac{1}{n}}L^{\frac{1}{n-1}}.\]
Let $m_{ij}(x)$ denote the order of vanishing of $q_{ij}$ at the point $x$. We know that 
\[\sum_j m_{ij}(x)\geq n_i(x)/2.\]
We choose non-negative numbers $n_{ij}(x)$ at each joint, satisfying
\begin{itemize}
\item $\sum n_{ij}(x)=n_i(x)/2$,
\item $n_{ij}(x)\leq m_{ij}(x)$.
\end{itemize}
For subsets $\sL'\subset \sL_i'$ and $\cJ'\subset \cJ_i'$ we define the weighted incidence count as
\[I_{ij}(\sL',\cJ')=\sum_{l\in \sL'}\sum_{x\in l\cap \cJ'} n_{ij}(x).\]
We know that
\[\sum_j I_{ij}(\sL_i',\cJ_i')=\sum_{x\in \cJ_i'}n_i(x)r_i(x)/2\geq d_iK^{\frac{n-1}{n}}L/16.\]
Since $\sum_j d_{ij}< d_i$, there exists a $j$ satisfying
\[I_{ij}(\sL_i',\cJ_i')\geq d_{ij}K^{\frac{n-1}{n}}L/16.\]
We denote the corresponding hypersurface by $S_2$. Similarly, we choose the sets $\sL_{ij}'$ etc. as we did in {\em Step} 2. We note that if $M$ is large enough then we can ensure that any point of $\cJ_{ij}$ and of $\cJ_{ij}'$ is a joint with respect to $\sL_{ij}$ and $\sL_{ij}'$ respectively.
\medskip

{\em Step} 4:
We keep going until we get $S_{n-2}$. Notice that the hypersurfaces $S_i$ are all irreducible and their degrees are strictly decreasing. Therefore they have no common component. We apply the ''Recursive step'' to the hypersurface $S_{n-2}$ as well, we obtain the subsets $\sL'\subset \sL$ and $\cJ'\subset \cJ$, and the quantities $q'$, $d'$, $n'(x)$, $r'(x)$. We remark that tracing through the recursive steps we obtain that $r'(x)>c_nM$ for all points $x\in \cJ'$ where $c_n$ is a constant depending only on $n$.

\noindent
Similarly to inequality \ref{eq:main} we have an inequality of the form
\[\sum_{x\in \cJ'}n'(x)r'(x)\geq c_nd'K^{\frac{n-1}{n}}L\]
where $c_n$ is a constant depending only on $n$. By H\"older's inequality we have
\[\left(\sum_{x\in \cJ'} n'(x)^n\right)^{\frac{1}{n}}\left(\sum_{x\in \cJ'} r'(x)^{\frac{n}{n-1}}\right)^{\frac{n-1}{n}}\geq c_nd'K^{\frac{n-1}{n}}L.\]
We estimate $\left(\sum_{x\in \cJ'} n'(x)^n\right)^{\frac{1}{n}}$ using the trivial bound:
\[\left(\sum_{x\in \cJ'} n'(x)^n\right)^{\frac{1}{n}}\leq \left(\sum_{x\in \cJ} r(x)^{\frac{n}{n-1}}\right)^{\frac{1}{n}}=2K^{\frac{1}{n}}L^{\frac{1}{n-1}}\]
obtaining
\[\left(\sum_{x\in \cJ'} r'(x)^{\frac{n}{n-1}}\right)^{\frac{n-1}{n}}\geq c_nd'K^{\frac{n-2}{n}}L^{\frac{n-2}{n-1}}\]
where $c_n$ is some constant depending only on $n$. Simplifying the inequality, we have
\begin{align}
\label{main}
\sum_{x\in \cJ'} r'(x)^{\frac{n}{n-1}}\geq c_n(d')^{\frac{n}{n-1}}K^{\frac{n-2}{n-1}}L^{\frac{(n-2)n}{(n-1)^2}}.\end{align}

We separate cases:
\begin{itemize}
\item $d'>nL^{\frac{1}{n-1}}$: In this case, we apply Lemma \ref{lem:lines}: there exists another surface $S$ with degree 
\[e<nL^{\frac{1}{n-1}}<d'\]
vanishing on all the lines in $\sL$. The hypersurface $S$ does not have a common component with the hypersurfaces $S_j$ ($1\leq j\leq n-2$), because the $S_j$ are irreducible and their degrees are at least $d'$ which satisfies $d'>e$. On the other hand since $r'(x)>c_nM$ for all points $x\in \cJ'$, we can apply Lemma \ref{lem:Kol} to the subset of points $\cJ'\subset S\cap S_1\cap...\cap S_{n-2}$ provided that $M$ is big enough:
\begin{align*}\sum_{x\in \cJ'} r'(x)^{\frac{n}{n-1}}&\leq ne\prod_{j=1}^{n-2} \deg(S_j)\cdot \max (\deg (S_j))\leq\\
&\leq n^2L^{\frac{1}{n-1}}d'd_i^{n-2}
\leq 2^{n-2}n^2d'K^{\frac{n-2}{n}}L.\end{align*}
A straightforward calculation shows that inequality \ref{main} would imply that 
\[d'\leq c_nK^{-\frac{n-2}{n}}L^{\frac{1}{n-1}}\]
where $c_n$ is some constant depending on only $n$. If we choose $K$ to be big enough, we get a contradiction to the assumption $d'>nL^{\frac{1}{n-1}}$.
\item $d'\leq nL^{\frac{1}{n-1}}$: We consider the gradient of $q'$ again. Choosing $M$ large enough we can ensure that all the points in $\cJ'$ are joints with respect to the set of lines $\sL'$, and thus all points in $\cJ'$ are singular points of $S_1\cap...\cap S_{n-2}$. Therefore, the components of the gradient vanish at those points. Since the degree of the components of the gradient is at most $d'-1$, there exists a component which does not vanish on all $S_1\cap...\cap S_{n-2}$. We choose the hypersurface $S$ to be the vanishing locus of such component. We apply Lemma \ref{lem:Kol} to the subset of points $\cJ'\subset S\cap S_1\cap...\cap S_{n-2}$ provided that $M$ is big enough as in the previous case:
\begin{align*}\sum_{x\in \cJ'} r'(x)^{\frac{n}{n-1}}&\leq n(d'-1)\prod_{j=1}^{n-2} \deg(S_j)\max(\deg(S_j))\leq\\
&\leq nd'^2d_i^{n-2}\leq 2^{n-2}n^2d'^2L^{\frac{n-2}{n-1}}K^{\frac{n-2}{n}}.\end{align*}
A straightforward calculation shows that the above inequality combined with Inequality \ref{main} implies that $d\geq c_nK^{\frac{1}{n}}L^{\frac{1}{n-1}}$ where $c_n$ is again some constant depending on only $n$. Choosing $K$ to be big enough, we obtain a contradiction to the assumption $d'\leq nL^{1}{n-1}$.
\end{itemize}
This concludes the proof of Theorem \ref{thm:main}.\qed
\end{Proof}

\section{Further remark}

 We remark that Theorem \ref{thm:main} fails without the extra hypothesis. The stronger inequality
\begin{align}
\label{ine:que}
\sum_{x\in \cJ}r(x)^{\frac{n}{n-1}}\leq c_nL^{\frac{n}{n-1}}
\end{align}
fails, even if $n=3$. We give two counterexamples:
\begin{itemize}
\item Consider a plane $P$ in $\F_p^3$. At each point of the plane $P$ take a line which is not parallel to $P$. Let $\sL$ be the set of these lines and all the lines contained in $P$. Clearly, $|\sL|$ is approximately $p^2$, hence the right-hand side of inequality \ref{ine:que} is approximately $p^3$. Lines in $\sL$ form joints at each point of $P$, and moreover at each point of $P$ there are approximately $p$ lines going through. Therefore, the left-hand side of inequality \ref{ine:que} is approximately $p^2\cdot p^{3/2}>p^3$.
\item Similar construction can be done using the Heisenberg surface $S\subset \F_{p^2}^3$ cut out  by the equation
\[x-x^p+y^pz-yz^p=0.\]
The Heisenberg surface contains $p^4$ lines and approximately $p^5$ points all incident to approximately $p$ lines. Similarly as before one can attach $p^4$ extra lines to the Heisenberg surface making all the points of the surface joints. Using the union of the set of lines contained in the Heisenberg surface and set of the attached lines, we see that the left-hand side of inequality \ref{ine:que} is approximately $p^5\cdot p^{3/2}$, the right-hand side, on the other hand, is approximately $p^6<p^5\cdot p^{3/2}$.
\end{itemize}
All the counterexamples we know involve ''planar'' surfaces (\cite{ElHa}). We believe that Theorem \ref{thm:main} remains true (in the $n=3$ case) assuming the strong Wolff axiom: that no $L^{1/2}$ lines lie in a planar surface. It would be interesting to find an explicit relationship between Conjecture \ref{con:main} and this stronger version of Theorem \ref{thm:main}.

\end{document}